\documentclass[a4paper,10pt]{article}

\usepackage[utf8]{inputenc}
\title{Finitude cohomologique des morphismes propres en géométrie algébrique : une preuve transcendante sans techniques projectives}
\author{\sc Antoine Ducros\\ \small Laboratoire J.-A. Dieudonné \\ \small Université de Nice- Sophia Antipolis, Parc Valrose\\ \small 06108 Nice CEDEX 02~~FRANCE}
\date{}
\usepackage{mathrsfs}

\usepackage[francais]{babel}
\usepackage{amssymb}
\usepackage{hyperref} 
\date{}
\input{xypic}

\newcommand{\hotimes}{\widehat{\otimes}}
\renewcommand{\phi}{\varphi}

\newcommand{\got}[1]{{\mathfrak #1}}
\renewcommand{\Bbb}{\mathbb}    
\renewcommand{\cal}{\mathscr}

\renewcommand{\H}{\mbox{\rm H}}
\newcommand{\RR}{{\Bbb R}}

\newcommand{\ZZ}{{\Bbb Z}}

\newcommand{\PP}{{\Bbb P}}

\newcommand{\grot}{_{\tiny \mbox{\rm G}}}

\newcommand{\QQ}{{\Bbb Q}}

\renewcommand{\epsilon}{\varepsilon}

\newcommand{\zero}{^{\mbox{\tiny o}}}

\newcommand{\cec} {\v{C}ech}
\newcommand{\itb}{\medskip \item[$\bullet$]}

\newcounter{cpt}
\newcounter{cptbis}

\begin{document}

\maketitle

\section*{\bf Introduction}

Si $f : {\cal Y} \to {\cal X}$ est un morphisme propre entre deux schémas localement noethériens, et si $\cal F$ est un ${\cal O}_{\cal Y}$-module cohérent, alors pour tout entier $q$ le ${\cal O}_{\cal X}$-module  $\mbox{R}^{q}f_{*}{\cal F}$ est cohérent. Pour démontrer ce résultat, on procède usuellement en deux temps ; on suppose tout d'abord que $f$ est projectif, et on utilise tout un arsenal spécifique à cette situation : tensorisation par ${\cal O}(n)$, calculs explicites des groupes de cohomologie de \cec ~relatifs au recouvrement affine standard de $\PP^{r}$, {\em etc.} ; on se ramène ensuite au cas projectif à l'aide du lemme de Chow. 

\medskip
On se propose dans ce qui suit de donner, lorsque $\cal Y$ et $\cal X$ sont localement de type fini sur un corps $k$, une nouvelle preuve de cette assertion qui n'utilise pas de techniques projectives et ne requiert pas le lemme de Chow. Le principe est très simple : 

\begin{itemize}

\itb Ce théorème est vrai dans le cadre des espaces analytiques au sens de Berkovich ; sa démonstration consiste à se ramener par changement du corps de base au cas d'un morphisme propre entre deux espaces analytiques rigide ; la finitude cohomologique d'un tel morphisme a été établie par Kiehl, sans recours au lemme de Chow bien entendu (un espace analytique rigide propre est loin, en général, d'être biméromorphiquement isomorphe à une variété projective) : Kiehl utilise des raffinements forts de recouvrements affinoïdes, et le caractère complètement continu des applications qu'un tel raffinement induit au niveau des complexes de \cec . 

\itb À toute variété $\cal Z$ sur $k$ peut être associée un espace analytique ${\cal Z}^{\beth}$ sur le corps $k$ {\em muni de la valeur absolue triviale}. Si $\cal {Y}\to \cal {X}$ est propre, alors $\cal {Y}^{\beth}\to \cal {X}^{\beth}$ est propre, d'après un résultat de Temkin qui repose sur des techniques valuatives (et ne fait pas appel au lemme de Chow.) 

\itb Un calcul explicite de complexes de \cec ~montre que le théorème de comparaison GAGA vaut dans ce cadre, ce qui permet de conclure. 
\end{itemize}

\section*{Démonstration du théorème}

\subsection*{Brefs rappels sur les faisceaux cohérents en géométrie analytique} 

Soit $k$ un corps ultramétrique complet, soit $Z$ un espace $k$-analytiques au sens de Berkovich (\cite{brk1}, \cite{brk2}). Soit $\cal F$ un faisceau cohérent sur le site $Z\grot$ (\cite{brk2}, \S 1.3). Si $Z$ est paracompact et analytiquement séparé ({\em i.e}, la diagonale est une immersion fermée), alors $Z$ possède un G-recouvrement affinoïde localement fini $\got{U}$ ; par le théorème d'acyclicité de Tate ({\em cf.} \cite{brk1}, prop. 2.2.5), et compte-tenu du fait que l'intersection de deux domaines affinoïdes de $Z$ est affinoïde en vertu de l'hypothèse de séparation, l'on dispose d'un isomorphisme naturel $\H^{\bullet}(Z\grot,{\cal F})\simeq \H^{\bullet}{\cal C}(\got{U},{\cal F})$, où $\cal C$ désigne le complexe de \v{C}ech. 

\medskip
Supposons que $Z$ est compact et $\got{U}$ fini. Pour tout domaine affinoïde $V$ de $Z$, munissons ${\cal F}(V)$ de sa structure de ${\cal A}_{V}$-module de Banach fini (\cite{brk1}), \S 2) ; ce choix définit pour tout $q$ une structure de $k$-espace vectoriel de Banach sur $\H^{q}{\cal C}(\got{U},{\cal F})$, et donc sur $\H^{q}(Z\grot,{\cal F})$. Si $L$ est une extension complète de $k$, le produit tensoriel complété par $L$ transforme suites exactes courtes admissibles en suites exactes courtes admissibles (Gruson, {\em cf.} \cite{brk2}, preuve du lemme 2.1.2) ; on en déduit l'existence pour tout $q$ d'un isomorphisme $\H^{q}((Z_{L})\grot,{\cal F}_{L})\simeq \H^{q}(Z\grot,{\cal F})\hotimes{L}$, où $Z_{L}$ et ${\cal F}_{L}$ ont le sens que l'on devine.

\medskip
Soit $f: Y\to X$ un morphisme propre, c'est-à-dire séparé, compact et sans bord, entre espaces $k$-analytiques et soit $\cal F$ un faisceau cohérent sur $Y\grot$. Soit $q$ un entier. Rappelons très brièvement comment Berkovich montre (\cite{brk1}, prop. 3.3.5), en s'appuyant sur un résultat de Kiehl, la cohérence de $\mbox{R}^{q}f_{*}{\cal F}$. On se ramène aussitôt au cas où $X$ est affinoïde, on écrit $X={\cal M}({\cal A})$ ; l'espace $Y$ est dès lors compact et il existe une famille ${\bf r}=(r_{1},\ldots,r_{n})$ de réels strictement positifs qui est libre dans $\QQ\otimes_{\ZZ}(\RR^{*}_{+}/|k^{*}|)$, et qui est telle que $|k_{\bf r}^{*}|\neq\{1\}$ et que $Y_{k_{\bf r}}$ et $X_{k_{\bf r}}$ soient {\em strictement} $k_{\bf r}$-analytiques (pour le sens de la notation 
$k_{\bf r}$, {\em cf.} \cite{brk1}, \S 2.1). Le morphisme $Y_{k_{\bf r}}\to X _{k_{\bf r}}$ est propre (\cite{brk2}, prop. 2.5.8 $ii)$). En vertu d'un théorème de Kiehl (\cite{kieeig}, th. 3.3),  $\mbox{R}^{q}f_{*}{\cal F}_{k_{\bf r}}$ est cohérent ; notons que la notion de propreté au sens de Berkovich correspond bien à celle définie par Kiehl, par exemple d'après le lemme 2.5.11 de \cite{brk1}. En conséquence : 

\medskip
\begin{itemize}
\item[$i)$] $\H^{q}((Y_{k_{\bf r}})\grot,{\cal F}_{k_{\bf r}})$ est un ${\cal A}\hotimes k_{\bf r}$-module de type fini, qui s'identifie par ailleurs comme on l'a vu à $\H^{q}(Y\grot,{\cal F})\hotimes k_{\bf r}$ ; on en déduit, à l'aide de la description explicite du corps $k_{\bf r}$ comme une algèbre de séries, que $\H^{q}(Y\grot,{\cal F})$ est un $\cal A$-module de type fini.

\item[$ii)$] Par un raisonnement analogue, on montre que si $V$ est un domaine affinoïde de $X$, alors 
\ $\H^{q}(f^{-1}(V)\grot,{\cal F})$ est isomorphe à $\H^{q}(Y\grot,{\cal F})\otimes {\cal A}_{V}$. 

\end{itemize}

\subsection*{L'espace analytique associé à un schéma et le théorème de comparaison} 

Soit $k$ un corps quelconque. À tout $k$-schéma de type fini ${\cal Z}$ on peut associer de deux manières différentes un espace analytique sur le corps $k$ {\em muni de la valeur absolue triviale} : la première consiste à considérer son analytification au sens de \cite{brk2}, \S 2.6 ; la seconde, qui a déjà été utilisée par Thuillier dans \cite{thuifan}, consiste à voir $\cal Z$ comme un schéma formel localement de présentation finie sur l'anneau discret $k$ (si l'on préfère, cela signifie que l'on s'intéresse à l'espace topologiquement annelé déduit de $\cal Z$ en munissant ${\cal O}_{\cal Z}(U)$ de la topologie discrète pour tout ouvert $U$ de $\cal Z$), et à considérer sa {\em fibre générique} au sens de \cite{vcy}. C'est ce que nous ferons ici, et nous noterons ${\cal Z}^{\beth}$ cette fibre générique. 

\medskip
Rappelons comment elle est cosntruite : si ${\cal Z}$ est affine d'anneau $\mathsf A$, alors ${\cal Z}^{\beth}$ est l'espace affinoïde ${\cal M}(\mathsf A^{\beth})$, où $\mathsf A^{\beth}$ est l'algèbre affinoïde égale à $\mathsf A$ munie de la norme triviale. Si $\mathsf A$ est une $k$-algèbre de type fini et si $f$ est un élément de $\mathsf A$, l'application $\mathsf A^{\beth} \to \mathsf A_{f}^{\beth}$ est évidemment bornée, et ${\cal M}(\mathsf A_{f}^{\beth})\to {\cal M}(\mathsf A^{\beth})$ identifie ${\cal M}(\mathsf A_{f}^{\beth})$ au domaine affinoïde d'équation $|f|=1$ de ${\cal M}(\mathsf A^{\beth})$. On utilise cette remarque pour construire ${\cal Z}^{\beth}$ dans le cas général par recollement à partir du cas affine ; on dispose d'une flèche de réduction ${\cal Z}^{\beth}\to {\cal Z}$,  {\em via} laquelle l'image réciproque d'un ouvert affine $V$ de $\cal Z$ s'identifie à $V^{\beth}$.

\medskip
Soit $\cal Z$ un schéma {\em séparé} et de type fini sur $k$. L'espace analytique ${\cal Z}^{\beth}$ est alors séparé et compact. Soit $(U_{i})$ un recouvrement fini de $\cal Z$ par des ouverts affines et soit $\cal F$ un faisceau cohérent sur $\cal Z$ ; pour tout $i$, on note $\mathsf A_{i}$ l'algèbre qui correspond à $U_{i}$, et $\mathsf M_{i}$ le $\mathsf A_{i}$-module de type fini ${\cal F}(U_{i})$. Chaque $\mathsf M_{i}$ peut être vu comme un $\mathsf A_{i}^{\beth}$-module de type fini. Les isomorphismes de compatibilité et d'associativité permettant de recoller les $\mathsf M_{i}$ perdurent une fois que l'on a muni les anneaux considérés de la norme triviale, ce qui permet de définir un faisceau cohérent ${\cal F}^{\beth}$ sur ${\cal Z}^{\beth}$ ; on vérifie immédiatement que ${\cal F}^{\beth}$ ne dépend pas du recouvrement ouvert choisi. 

\medskip
Par construction, ${\cal F}^{\beth}(V^{\beth})={\cal F}(V)$ pour tout ouvert affine $V$ de $\cal Z$. En particulier, les complexes de \cec ~${\cal C}((U_{i}),{\cal F})$ et ${\cal C}((U_{i}^{\beth}),{\cal F}^{\beth})$ sont isomorphes, d'où l'on déduit l'existence d'un isomorphisme $$\H^{\bullet}({\cal Z},{\cal F})\simeq \H^{\bullet}({\cal Z}^{\beth},{\cal F}^{\beth}).$$ 

\subsection*{Preuve du théorème}

\noindent
{\bf Théorème.} {\em Soit $k$ un corps et soit $f: {\cal Y}\to {\cal X}$ un morphisme propre entre $k$-schémas localement de type fini. Soit $\cal F$ un faisceau cohérent sur $\cal Y$. Pour tout entier $q$, le ${\cal O}_{\cal X}$-module $\mbox{\rm R}^{q}f_{*}{\cal F}$ est cohérent.} 

\bigskip
{\em Démonstration.} On peut supposer $\cal X$ affine ; soit $\mathsf A$ l'anneau correspondant. Comme  $\mbox{R}^{q}f_{*}{\cal F}$ est quasi-cohérent, il suffit de vérifier que le $\mathsf A$-module $\H^{q}({\cal Y},{\cal F})$ est de type fini. Il résulte du paragraphe précédent que $\H^{q}({\cal Y},{\cal F})$ est isomorphe à $\H^{q}({\cal Y}^{\beth},{\cal F}^{\beth})$ ; comme l'anneau ${\mathsf A}^{\beth}$ coïncide avec $\mathsf A$, il suffit de montrer que $\H^{q}({\cal Y}^{\beth},{\cal F}^{\beth})$ est un ${\mathsf A}^{\beth}$-module de type fini. On va montrer que ${\cal Y}^{\beth}\to {\cal X}^{\beth}$ est propre, ce qui permettra de conclure.

\medskip
Comme ${\cal Y}^{\beth}\to {\cal X}^{\beth}$ est automatiquement compact, il suffit de montrer qu'il est sans bord. Ceci peut se vérifier après extension des scalaires à $k_{r}$ pour un certain $r$ de $]0;1[$, par le corollaire 2.5.12 de \cite{brk1}. 

\medskip
Remarquons que ${\cal Y}^{\beth}_{k_{r}}$ (resp. ${\cal X}^{\beth}_{k_{r}}$) apparaît comme la fibre générique du schéma formel ${\cal Y}\times \mbox{Spf}\;k_{r}\zero$ (resp. ${\cal X}\times \mbox{Spf}\;k_{r}\zero$). Comme ${\cal Y}\to {\cal X}$ est propre, 
${\cal Y}\times \mbox{Spf}\;k_{r}\zero \to {\cal X}\times \mbox{Spf}\;k_{r}\zero$ est propre, puisque cela se teste par définition au niveau des fibres spéciales ; la valeur absolue de $k_{r}$ n'étant pas triviale, ${\cal Y}^{\beth}_{k_{r}}\to  {\cal X}^{\beth}_{k_{r}}$ est propre d'après un résultat de Temkin (\cite{tmk}, cor. 4.4), dont la preuve est purement valuative.~$\Box$

\end{document}